\documentclass[11pt]{article} 

\usepackage{tikz}
\usetikzlibrary{positioning}
\usepackage[utf8]{inputenc} 
\usepackage{amsmath,breqn}
\usepackage{amssymb}
\usepackage{graphicx}
\usepackage{epsfig,epstopdf}
\usepackage{pdflscape}
\usepackage{rotating,lscape,afterpage}
\usepackage{hyperref}

\usepackage{geometry} 
\usepackage{graphicx} 

\usepackage{booktabs} 
\usepackage{array} 
\usepackage{paralist} 
\usepackage{verbatim} 
\usepackage{subfig} 

\usepackage{fancyhdr} 
\usepackage{sectsty}
\allsectionsfont{\sffamily\mdseries\upshape} 

\usepackage[nottoc,notlof,notlot]{tocbibind} 
\usepackage[titles,subfigure]{tocloft} 

\title{Equilibrium solutions of three player Kuhn poker with $N>3$ cards: \\A new numerical method using regularization and arc-length continuation}
\author{John Billingham\\ School of Mathematical Sciences, \\ The University of Nottingham, \\ Nottingham NG7 2RD, UK}

\begin{document}
\begin{abstract}
We study the equilibrium solutions of three player Kuhn poker with $N>3$ cards. We compute these solutions as a function of the initial pot size, $P$, using a novel method based on regularizing the system of polynomial equations and inequalities that defines the solutions, and solving the resulting system of nonlinear, algebraic equations using a combination of Newton's method and arc-length continuation. We find that the structure of the equilibrium solution curve is very complex, even for games with a small number of cards. Standard three player Kuhn poker, which is played with $N=4$ cards, is qualitatively different from the game with $N>4$ cards because of the simplicity of the structure of the value betting and bluffing ranges of each player. When $N>5$, we find that there is a new type of equilibrium bet with midrange cards that acts as a bluff against one player and a value bet against the other. 
\end{abstract}
\maketitle

\section{Introduction}\label{sec_intro}
The computation of Nash equilibria in multiplayer games is a central problem in computational game theory (see, for example, \cite{vonS2010} and the papers described therein). Many algorithms exist, with different strengths and weaknesses, suitable for a variety of different forms of game (for example, see the software package {\em Gambit}, \cite{Gambit}). In this paper, we study a simple class of toy poker games - three player Kuhn poker with $N>3$ cards and pot size $P$ - and focus on computing the equilibrium solutions as $N$ and $P$ vary. The structure of the system of polynomial equations and inequalities that determines the equilibrium solutions is very simple, and we have developed an algorithm that exploits this.

Three player Kuhn poker has previously been studied with $N=4$ and $P=3$, which is one of the simplest multiplayer toy poker games that can be used to test the performance of equilibrium finding algorithms and more general game playing agents \cite{Risk:2010:UCR:1838206.1838229,Szafron:2013:PFE:2484920.2484962}. We show in this paper that the structure of the equilibrium solution curve as a function of $P$ is complex, and that for $N=4$ the initial pot size $P=3$ is a special case for which a range of equilibrium solutions is possible. We also find that increasing $N$ to five and beyond leads to a significant increase in the complexity of the equilibrium solution curve.

It is worth making some comments about the philosophy behind our investigation. We believe that toy poker games such as those examined in this paper are interesting in their own right. They are mathematical objects of surprising complexity, given their simple definition. There is also value in studying and understanding the simplest nontrivial multiplayer games as a basis for understanding more complex multiplayer games with greater significance to the real world. In addition, poker (Texas Holdem in particular) is one of the world's most popular games, and is often cited as a challenge problem in game theory and computer science (for example, \cite{Billings:2002:CP:512148.512158}). Although an equilibrium solution of two player, limit Texas Holdem has recently been computed, \cite{15science}, and an AI based on a combination of equilibrium finding and neural network evaluation performs at a superhuman level in two player, no limit Texas Holdem, \cite{DSAI}, the multiplayer versions of these games have not been studied in as much detail, \cite{Risk:2010:UCR:1838206.1838229}. 

In Section~\ref{sec_game}, we describe the class of game, three player Kuhn poker with $N>3$ cards, that we will study in this paper (most easily understood with reference to the decision tree shown in Figure~\ref{fig_decision}), and formulate the mathematical problem that determines the equilibrium solutions. We also note the dominated strategies that can be eliminated, and determine the smallest pot size for which a nontrivial solution exists. In Section~\ref{sec_reg}, we introduce a regularized version of the equilibrium equations, and describe our implementation of the arc-length continuation method that we use to solve it. After verifying our method on a related game that has an analytical solution, namely simplified Kuhn poker (SKP, \cite{SKP2017_2}), the results of which are discussed in Section~\ref{sec_SKP}, we move on to study numerical solutions of the full problem for various $N>3$ in Section~\ref{sec_num}. We discuss our results further and indicate some extensions of this work that we will pursue in the future in Section~\ref{sec_conc}.

\section{Three player Kuhn poker with $N>3$ cards}\label{sec_game}
In the version of three player Kuhn poker that we study in this paper, the deck contains $N>3$ different cards numbered from $1$ to $N$. A single card from the deck is dealt at random to each player without replacement. The pot contains $P$ units (each player is assumed to have contributed $P/3$ units to the pot). The possible betting sequences are shown in Figure~\ref{fig_decision}. The first action is made by Player 1, who can either bet one unit or check. The next action is made by Player 2, who can either call the bet of one unit or fold after Player 1 bets, and either check or bet one unit after Player 1 checks. The action then moves to Player 3, and continues in the same manner. The game ends after either three checks or after all players have either bet, called or folded. If two players fold, the remaining player wins the pot of $P$ units. If two or three players do not fold, the remaining player who holds the card with the largest numerical value wins the pot of $P$ units along with any bets or calls made by themselves and the other players. 
\afterpage{\begin{landscape}
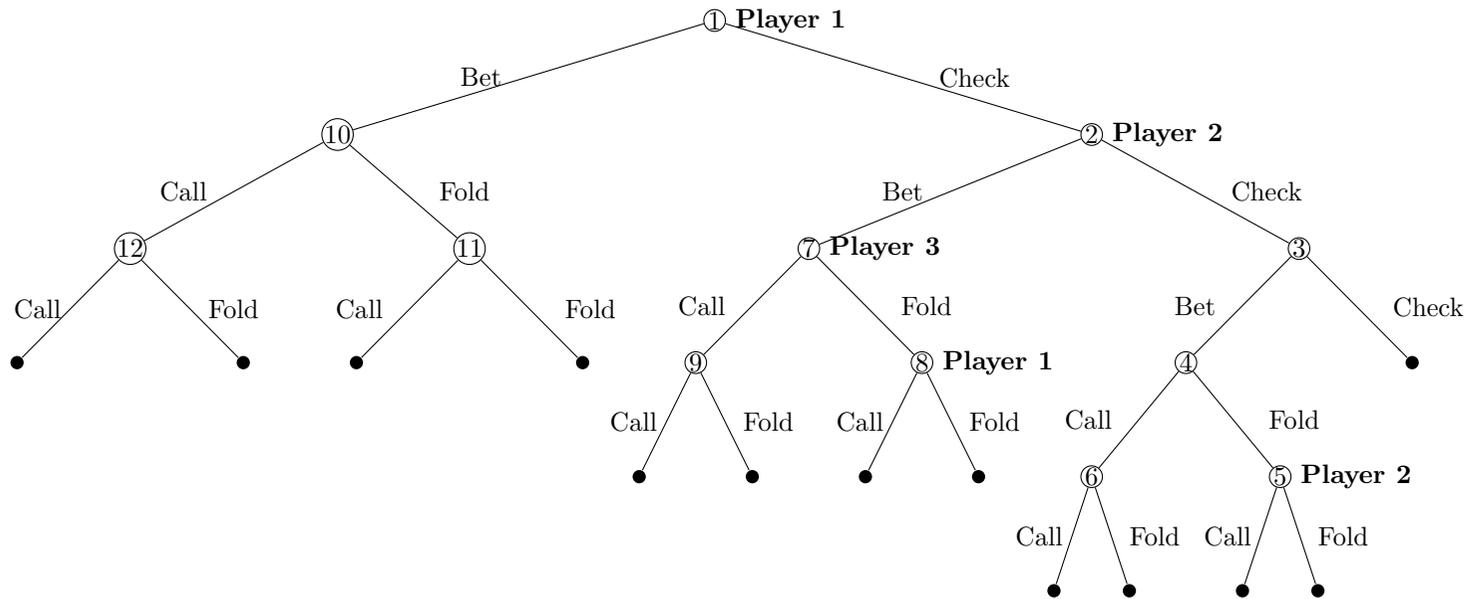
\begin{figure}
 	\begin{center}
    \small
    \begin{tikzpicture}[thin,
      level 1/.style={sibling distance=100mm},
      level 2/.style={sibling distance=55mm},
      level 3/.style={sibling distance=30mm},
      level 4/.style={sibling distance=15mm},
      level 5/.style={sibling distance=10mm},
      every circle node/.style={minimum size=1.75mm,inner sep=0mm}]

\node[circle,draw,label=right:{\bf Player 1}] {1}
child {
	node [circle,draw] {10} 
		child {
			node[circle,draw] {12}
			child {
				node [circle,fill] {} edge from parent node[left] {Call~}
				}
			child {
				node [circle,fill] {} edge from parent node[right] {~Fold}
				} edge from parent node[left] {Call~~~}
			}
		child {
			node[circle,xshift = -10mm,draw] {11}
			child {
				node [circle,fill] {} edge from parent node[left] {Call~~~}
				}
			child {
				node [circle,fill] {} edge from parent node[right] {~~~Fold}
				} edge from parent node[right] {~~~Fold}
			}edge from parent node[left] {Bet~~~}
	}
child {
	node [circle,draw,label=right:{\bf Player 2}] {2}
		child {
			node[circle,xshift = -10mm,draw,label=right:{\bf Player 3}] {7}
			child {
				node [circle,draw] {9}
			child {
				node [circle,fill] {} edge from parent node[left] {Call~}
				}
			child {
				node [circle,fill] {} edge from parent node[right] {~Fold}
				} edge from parent node[left] {Call~~~}
				}
			child {
				node [circle,draw,label=right:{\bf Player 1}] {8}
			child {
				node [circle,fill] {} edge from parent node[left] {Call~}
				}
			child {
				node [circle,fill] {} edge from parent node[right] {~Fold}
				} edge from parent node[right] {~~~Fold}
				} edge from parent node[left] {Bet~~~}
			}
		child {
			node[circle,draw] {3}
			child {
				node [circle,draw] {4}
		child {
			node[circle,draw,xshift=-5mm] {6}
			child {
				node [circle,xshift=0mm,fill] {} edge from parent node[left] {Call~}
				}
			child {
				node [circle,xshift=0mm,fill] {} edge from parent node[right] {~Fold}
				} edge from parent node[left] {Call~~~}
			}
		child {
			node[circle,draw,label=right:{\bf Player 2},xshift=5mm] {5}
			child {
				node [circle,xshift=0mm,fill] {} edge from parent node[left] {Call~}
				}
			child {
				node [circle,xshift=0mm,fill] {} edge from parent node[right] {~Fold}
				} edge from parent node[right] {~~~Fold}
			} edge from parent node[left] {Bet~~~}
				}
			child {
				node [circle,fill] {} edge from parent node[right] {~~~Check}
				} edge from parent node[right] {~~~Check}
			} edge from parent node[right] {~~~Check}
	}
;

    \end{tikzpicture}
    \end{center}
    \caption{The decision tree for three player Kuhn poker. Open circles are decision nodes, numbered from 1 to 12. Solid circles are terminal nodes. The Player labels given in bold face apply to all nodes on the same horizontal level, i.e. Player 1 acts at nodes 1, 4, 8 and 9, Player 2 acts at nodes 2, 5, 6 and 10, Player 3 acts at nodes 3, 7, 11 and 12.\label{fig_decision}}
  \end{figure}
\end{landscape}}

In this game, each player can choose between exactly two options (bet or check, call or fold) at each decision node at which they act. We denote by $x_{ij}$ the frequency at which an aggressive action (bet or call) is taken by Player $i$ for $i = 1$, $2$, $3$ at the four nodes that they control (see Figure~\ref{fig_decision}) with each of the $N$ possible cards, so that $j = 1$, $2$, $\ldots$ $4N$. The frequencies at which Player $i$ chooses a passive action (check or fold) are therefore $1-x_{ij}$. We will also use the notation $x_{ij} = x_l$ with $l = 4N(i-1)+j = 1$, $2$, $\ldots 12N$ whenever convenient for the presentation of our results.

The value of the game to Player $k$ is $E_k(x_{ij})$ for $k = 1$, $2$ or $3$. Note that this is a zero sum game, so $E_1+E_2+E_3 = 0$. An analytical expression for these three functions can be calculated, either by hand or using computer algebra, in terms of the frequencies $x_{ij}$, or calculated numerically from the game tree for a given set of $12N$ numerical values of $x_{ij}$. In this paper, we have used computer algebra (specifically, the Symbolic Math Toolbox in MATLAB) to determine $E_k$ and its derivatives analytically, and will discuss the limitations of this approach below. Note that each of the polynomial functions $E_k$ is linear in each of the frequencies $x_{ij}$. 

Since Player $k$ seeks to maximise $E_k$, an equilibrium set of betting frequencies, $x_{ij} = \hat{x}_{ij}$ satisfies
\begin{equation}
\frac{\partial E_i}{\partial x_{ij}}\left(\hat{x}_{ij}\right) = 0,~~\mbox{or}~~ \frac{\partial E_i}{\partial x_{ij}}\left(\hat{x}_{ij}\right) <0~~\mbox{and}~~\hat{x}_{ij}=0,~~\mbox{or}~~\frac{\partial E_i}{\partial x_{ij}}\left(\hat{x}_{ij}\right) >0~~\mbox{and}~~\hat{x}_{ij}=1~~\mbox{for $j = 1$, $2$, $\ldots 4N$}.\label{eqn1}
\end{equation}
This set of equations and inequalities (which can be rewritten as a variational inequality, although we will not do so here) takes this simple form because the decision tree has just two branches at each node. It is (\ref{eqn1}) that we study in this paper. We seek equilibrium solutions and wish to determine how they vary with pot size, $P$. 

In \cite{SKP2017_2}, we studied simplified Kuhn poker (SKP), in which there are just eleven non-trivial betting frequencies, and hence eleven conditions in (\ref{eqn1}). The equilibria of SKP can be found by using computer algebra to consider each of the $3^{11} = 177147$ possible combinations of conditions in (\ref{eqn1}). As the number of conditions in (\ref{eqn1}) increases, analysis of the equilibrium solutions rapidly becomes intractible using this exhaustive method. In order to study larger games, such as three player Kuhn poker ($N=4$), which has a system of 48 conditions in its raw state that can be reduced to 23 conditions by eliminating dominated strategies (see \cite{Szafron:2013:PFE:2484920.2484962}), and also the game for $N>4$, a numerical method must be used.  We propose a regularized version of (\ref{eqn1}) that transforms the problem into a system of nonlinear equations, which can be solved using Newton's method and arc-length continuation (see Section~\ref{sec_reg}).

Before proceeding, note that there are some dominated strategies that we can remove.
\begin{itemize}
\item A player holding the best card, with value $N$, will always bet or call at nodes 3 to 12. 
\item A player holding the worst card, with value $1$, will always fold at nodes 4 to 12. 
\item A player holding the second worst card, with value $2$, will always fold at nodes 6, 9 and 12.
\end{itemize}
This fixes twenty two of the betting frequencies, which leaves $12N-22$ equilibrium betting frequencies to be determined by (\ref{eqn1}). When $N=4$, it is also straightforward to show that checking with the second best card, with value $3$, dominates betting at nodes 1, 2 and 3. We will not enforce this in the results shown below.

Whatever the value of $N$, we find that an equilibrium strategy consists of bets with some high value cards and some low value cards at nodes 1, 2 and 3 (value bets and bluffs) and calls with some mid to high value cards at the other nodes (bluff catching). In addition, some fraction of the high value cards is checked at nodes 1 and 2 in order to call with them at nodes further down the decision tree (sandbagging). It is straightforward to show that when $N=4$, sandbagging and value betting occurs only with the highest value card, $4$, and bluff catching only with the second highest value card, $3$, but that bluffing may be with either of the two lowest value cards, $1$ or $2$, depending on the size of the pot, $P$. We will see below that, as the number of cards in the deck, $N$, increases, the structure of these ranges (value betting, bluffing, sandbagging, bluff catching) becomes more complex, and the multiplicity of the equilibrium solutions can become as high as twenty four for some values of $P$ (see Figure~\ref{E_N9}). We also find that there is a new type of bet, which acts as a bluff against one player and a value bet against the other, that obviously cannot exist in two player games.

Finally, we can use a simple argument to determine the smallest value of the pot size, $P$, for which any player can make a profit at equilibrium. There are $(N-1)(N-2)$ deals for which a given player holds the lowest value card and $(N-2)(N-3)$ deals for which a given player holds the lowest value card and the highest value card is not dealt. Bluffing with the lowest value card is potentially profitable if the money won when no strongest card is dealt and every player folds is greater than the money lost when a player with the strongest card calls, which gives $(N-2)(N-3)P > (N-1)(N-2)-(N-2)(N-3)$, and hence $P>P_{min}$, where
\[P_{min} \equiv \frac{2}{N-3}.\]
Our numerical results confirm that $P_{min}$ is indeed the pot size below which no profit is possible; more specifically, for all equilibrium solutions with $P \leq P_{min}$, $E_1 = E_2 = E_3 = 0$, and with $P>P_{min}$, $E_3>0$.

\section{The regularized game and numerical method} \label{sec_reg}
\begin{figure}
 \begin{center}
 \includegraphics[width=0.8\textwidth]{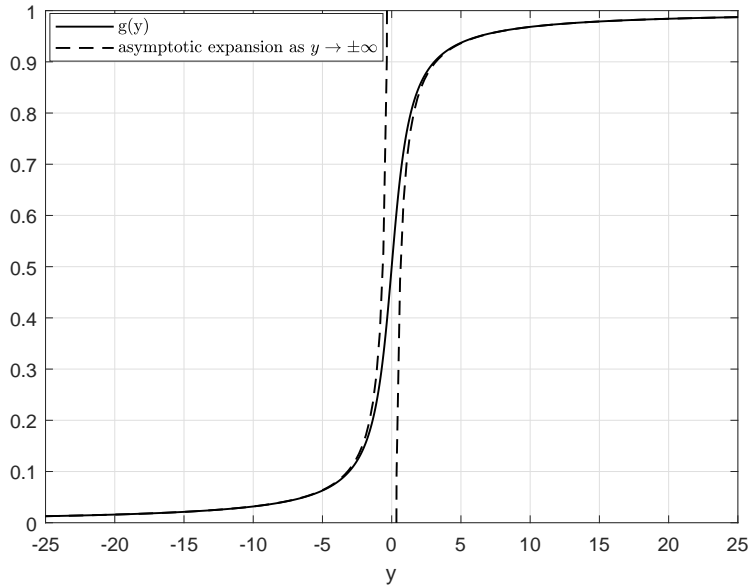}
 \caption{The regularization function, $g(y) = \frac{1}{2}+\frac{1}{\pi}\tan^{-1}y$, used in this paper, along with the asymptotic expansions $g \sim -1/\pi y$ as $y\to -\infty$ and $1-g \sim 1/\pi y$ as $y \to \infty$.}\label{fig_g}
 \end{center}
 \end{figure}

Consider the system of nonlinear equations 
\begin{equation}
g\left(\frac{1}{\epsilon}\frac{\partial E_i}{\partial x_{ij}}\right) - x_{ij} = 0,\label{eqn2}
\end{equation}
where $\epsilon$ is a positive regularization parameter and $g$ is a differentiable, monotonically increasing function that has
\begin{equation}
\begin{array}{c}
1- g(y) \sim k_+/y~~\mbox{as $y \to \infty$,} \\
g(y) \sim -k_-/y~~\mbox{as $y \to -\infty$,}
\end{array}\label{g_asym}
\end{equation}
with $k_{\pm}$ positive constants. This asymptotic behaviour is crucial for our method to work, for reasons explained in Appendix~\ref{sec_asol}. In all of the calculations discussed in this paper, we have used
\begin{equation}
g(y) = \frac{1}{2} + \frac{1}{\pi} \tan^{-1}(y),\label{eqn_g}
\end{equation}
which has $k_{\pm} = 1/\pi$ (as shown in Figure~\ref{fig_g}).
As $\epsilon \to 0$, (\ref{eqn2}) asymptotes to (\ref{eqn1}) in a manner that we investigate in more detail in Appendix~\ref{sec_asol}. 

Using the single suffix notation for the solution, it is convenient to rewrite (\ref{eqn1}) as 
\begin{equation}
f_k({\bf x})= 0,~~\mbox{or}~~ f_k({\bf x})<0~~\mbox{and}~~x_{k}=0,~~\mbox{or}~~f_k({\bf x}) >0~~\mbox{and}~~x_{k}=1,\label{eqn1_2}
\end{equation}
and hence the regularized system as
\begin{equation}
g\left(\frac{f_k({\bf x})}{\epsilon}\right) - x_{k} = 0,~~\mbox{for $k = 1$, $2 \ldots 12N-22$.}\label{eqn2_2}
\end{equation}
In order to solve (\ref{eqn2_2}) using Newton's method and arc-length continuation, we treat the pot size $P$ as an additional unknown and solve for the vector
\[{\bf X} \equiv \left({\bf x}, P\right).\]
Starting from two initial solutions, ${\bf X} = {\bf X}_0$ and ${\bf X} = {\bf X}_1$, we solve for a sequence of solutions ${\bf X}_i$, for $i = 2$, $3 \ldots$, augmenting the system (\ref{eqn2_2}) with the condition
\begin{equation}
\left({\bf X}_{i+1}-{\bf X}_i\right) \cdot \left( \frac{{\bf X}_{i}-{\bf X}_{i-1}}{\left\|{\bf X}_{i}-{\bf X}_{i-1}\right\|} \right)= \delta_{i+1},~~\mbox{for $i = 1$, $2 \ldots$,} \label{cont_eqn}
\end{equation}
where $\delta_i$ is a sequence of step sizes whose selection we discuss below. Equation (\ref{cont_eqn}) ensures that the size of the change in the solution vector from the solution at the previous step, ${\bf X}_{i+1} - {\bf X}_{i}$, projected in the direction tangent to the curve of solutions, is $\delta_{i+1}$. At each step, we solve (\ref{eqn2_2}) and (\ref{cont_eqn}) using Newton's method, with an  initial guess extrapolated quadratically from the previous three solutions. We precalculate the Jacobian analytically using the Symbolic Math Toolbox in MATLAB.

By regularizing the problem using the small parameter $\epsilon$, we transform it into a system of nonlinear equations. Although solving such systems is, in principle, more straightforward than solving systems of variational inequalities, we know that the solution when $\epsilon = 0$ is not smooth, and expect that regularization will lead to solutions smoothed over a scale determined by the size of $\epsilon$. This means that, during the arc-length continuation described above, there will be small regions where the solution changes rapidly. The sequence of step sizes, $\delta_i$, must therefore be chosen adaptively to take this into account. This is done by reducing $\delta_{i+1}$ and recomputing if either Newton's method fails to converge or $\left\|{\bf X}_{i+1}-{\bf X}_{i}\right\|$ is larger than $\mbox{max}\left\{1.05\delta_{i+1}, 10^{-3} \epsilon\right\}$, otherwise slightly increasing $\delta_{i+1}$, but not beyond a maximum value of $0.1$. In this manner, the sequence $\delta_i$ adapts to the contours of the solution curve. 

In order to find the initial solutions ${\bf X}_0$ and ${\bf X}_1$, we used the MATLAB solver \texttt{fsolve}, an implementation of the trust-region dogleg algorithm, to solve (\ref{eqn2_2}) with $P=0$ and $\epsilon = 0.1$, taking an initial guess for ${\bf x}$ drawn at random from a uniform distribution over $[0,1]$. This algorithm is able to find the solution using a poor initial guess in a manner that is not possible using Newton's method. We then use parameter continuation to reduce $\epsilon$ to the required value. We find the second initial solution using \texttt{fsolve} with $P=0.01$, taking the solution at $P=0$ as the initial guess. From these two initial solutions, we can proceed with arc-length continuation using Newton's method. 

Finally, for betting frequencies controlled by the same player with the same card that also lie in the same branch of the decision tree, for example, the betting frequencies of Player 1 at Nodes $1$ and $4$, we differentiated with respect to {\em both} betting frequencies in (\ref{eqn1}). This eliminates some equilibrium solutions that are not perfect in the subgame associated with the node lower down the tree. For example, if Player 1 always bets with card $1$ at Node $1$, $x_{11} = 1$. However, the expression for $\partial E_1/\partial x_{14}$ associated with the frequency of betting with card $1$ at node $4$ has a factor of $1-x_{11}$, and is therefore zero if $x_{11} = 1$, since Player 1 never reaches Node 4 with card $1$ in this case. The equilibrium will therefore not necessarily be correct in this case for the subgame at Node 4. By also differentiating with respect to $1-x_{11}$, i.e. using $-\partial^2 E_1/\partial x_{11}\partial x_{14}$, we remove the factor of $1-x_{11}$, and the equilibrium is then valid in the subgame at node 4. In our numerical investigation, we found that if we neither removed the dominated strategies listed earlier, nor performed these extra differentiations, spurious solutions, which use non-credible threats, could be computed. If either the dominated strategies are removed, or the extra differentiations are performed, these spurious solutions are eliminated. In the numerical results discussed below, we both removed the dominated strategies {\em and} performed the extra differentiation, since both steps are well-justified, but using either alone produces essentially the same results.

\section{Test Case: Simplified three player Kuhn poker}\label{sec_SKP}
In order to illustrate how solutions of the regularized system (\ref{eqn2}) converge to those of the original system (\ref{eqn1}) as $\epsilon \to 0 $, and also to verify our method against a known, analytical solution, we consider SKP, as studied in \cite{SKP2017_2}. In order to do this within the framework of three player Kuhn poker with $N$ cards, we set $N=4$ and force the betting and calling frequencies with the card of value $1$ to be zero. Note that the cards referred to as A, K, Q and J in \cite{SKP2017_2} correspond to $4$, $3$, $2$ and $1$ in the notation that we use here. 

\begin{figure}
 \begin{center}
 \includegraphics[width=0.85\textwidth]{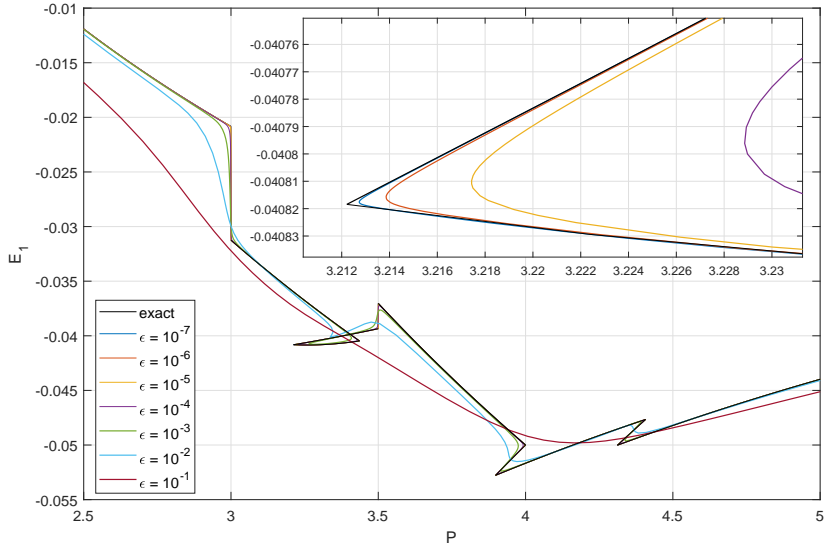}
 \caption{The exact expectation of Player 1 in SKP, along with the computed expectation for various values of $\epsilon$.}\label{fig_E1}
 \end{center}
 \end{figure}

\begin{figure}
 \begin{center}
 \includegraphics[width=0.85\textwidth]{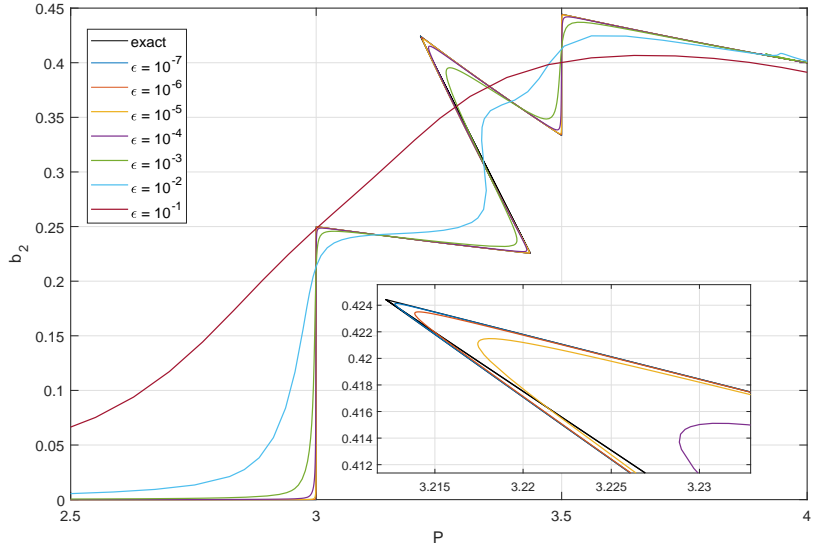}
 \caption{The exact equilibrium  betting frequency $b_2$, the bluffing frequency of Player 2 with a $2$ (Q) at Node 2, in SKP, along with the values computed for various values of $\epsilon$.}\label{fig_b2}
 \end{center}
 \end{figure}

The expectation of Player 1 is shown in Figure~\ref{fig_E1} for various values of the regularization parameter, $\epsilon$, along with the exact solution. Firstly, we can see that the use of arc-length continuation has allowed us to compute all of the equilibrium solutions, including the multiple equilibrium solutions that exist for some ranges of values of the pot size, $P$. Secondly, note that this expectation, which is an excellent indicator of the structure of the equilibrium solutions, showing clearly where multiple equilibrium solutions exist, is accurately reproduced by the regularized system as $\epsilon \to 0$. On the main scale shown, the numerical solution is indistinguishable from the analytical solution for $\epsilon < 10^{-4}$. If we zoom into a region where the exact expectation is not smooth, we can see how the numerically-calculated expectation becomes progressively more acccurate as $\epsilon$ decreases. Finally, we can see that at $P=3$ and $P=3.5$, where we know from the analytical solution that a range of different expectations is possible, regularization allows the expectations smoothly to connect in a single curve. Figure~\ref{fig_b2} shows that similar observations can be made about a typical betting frequency. 

\begin{figure}
 \begin{center}
 \includegraphics[width=0.8\textwidth]{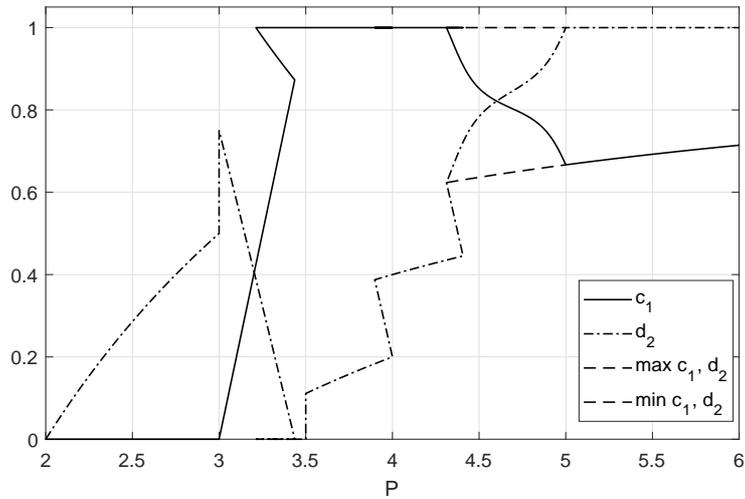}
 \caption{Computed betting frequencies $c_1$ and $d_2$ when $\epsilon = 10^{-6}$ for SKP. The exact upper and lower bounds on these frequencies are also shown.}\label{fig_c1d2}
 \end{center}
 \end{figure}

Figure~\ref{fig_c1d2} shows both the analytical equilibrium betting frequencies $c_1$ and $d_2$, which correspond to Player 1's calling frequency at Node 4 and Player 2's calling frequency at Node 5 with a card of value $3$ (a K), and those determined numerically with $\epsilon = 10^{-6}$. Also shown are the analytically-determined upper and lower bounds on $c_1$ and $d_2$ for the range of values of $P$ for which only $c_1+d_2$ is determined uniquely at equilibrium. As we can see, the regularized version of the problem selects unique values of $c_1$ and $d_2$ (see Appendix~\ref{sec_asol} for more discussion of this). 

\section{Numerical Results}\label{sec_num}
In this section, we will use the regularized problem with $\epsilon = 10^{-6}$ to investigate the equilibrium solutions of three player Kuhn poker with $N>3$ cards. In each case, we find a connected branch of equilibrium solutions. We have not attempted to find solutions not connected to this branch, and, although it seems plausible that no other solutions exist (this is the case for SKP), we have not attempted to investigate this further, preferring to focus on the structure of the branch of equilibrium solutions that we have computed.

\begin{figure}
 \begin{center}
 \includegraphics[width=\textwidth]{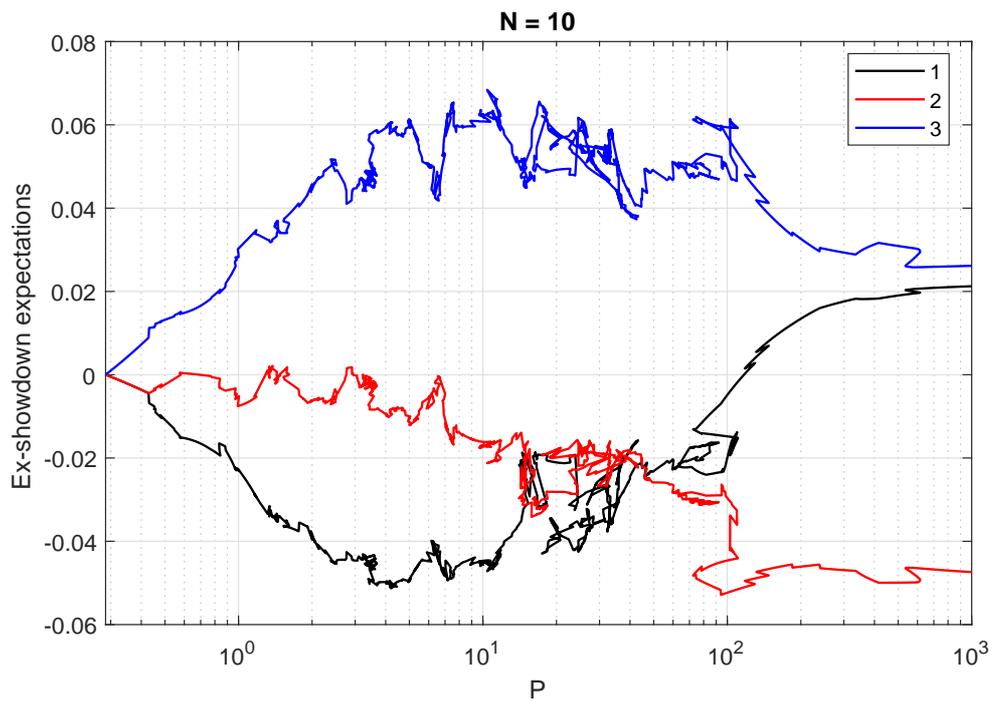}
 \caption{Computed expectations for $N=10$. Note that the scale for the pot size, $P$, is logarithmic.}\label{E_N10}
 \end{center}
 \end{figure}

We calculated the equilibrium solutions for many different values of $N$ with $P$ as large as $1000$. For example, Figure~\ref{E_N10} shows the expectations of each player when there are ten cards in the deck ($N=10$). This is a good way to visualize the solution structure, since different expectations in this plot must correspond to distinct equilibrium solutions. 
The most striking feature is the complexity of this solution curve, with up to fifteen distinct equilibrium solutions coexisting for some values of the pot size, $P$. Another surprising feature is that, for sufficiently large $P$, Player 2 becomes the biggest loser at equilibrium, with both $E_3$ and $E_1$ positive. For lower values of $P$, $E_3 > E_2 > E_1$, so the later the player acts (poker players say, the better the position of the player), the greater their profit at equilibrium. This switch to $E_2>E_3$ for $P$ greater than about $45$ is hard to explain and is a feature that exists when there are five or more cards (but not in standard three player Kuhn poker, which is played with four cards). 

Although the structure of the equilibrium solutions is clearly very complex for larger values of $P$, these games have little relevance to poker as played in the real world. Even with very aggressive betting on each street in games such as limit Texas Holdem, Seven Card Stud, Razz, Triple Draw and others, the pot to bet ratio, $P$, rarely reaches double figures. We therefore concentrate here on interesting features of the equilibrium solutions of three player Kuhn poker with $N>3$ and $P_{min}\leq P \leq10$. 

\subsection{Three player Kuhn poker ($N=4$)}
\begin{figure}
 \begin{center}
 \includegraphics[width=\textwidth]{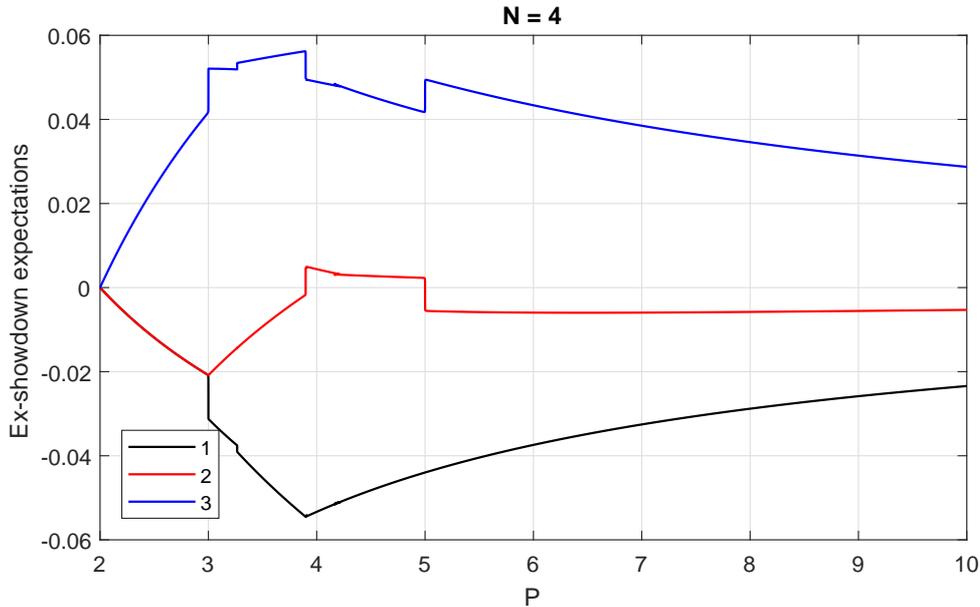}
 \caption{Computed expectations for $N=4$.}\label{E_N4_2}
 \end{center}
 \end{figure}

Figure~\ref{E_N4_2} shows the computed expectations for three player Kuhn poker ($N=4$). As we shall see, compared to games with more cards, the equilibrium solution structure is quite simple; even simpler than the solution structure of SKP, \cite{SKP2017_2}. There are four values of $P$, clearly visible in Figure~\ref{E_N4_2}, where a range of expectations is possible for two of the three players, including the case $P=3$, which is the usual pot size for three player Kuhn poker (each player antes one unit before play begins), and which was studied in \cite{Szafron:2013:PFE:2484920.2484962}. The analytical solutions found in \cite{Szafron:2013:PFE:2484920.2484962} are consistent with the solutions that we have computed. 

\begin{figure}
 \begin{center}
 \includegraphics[width=\textwidth]{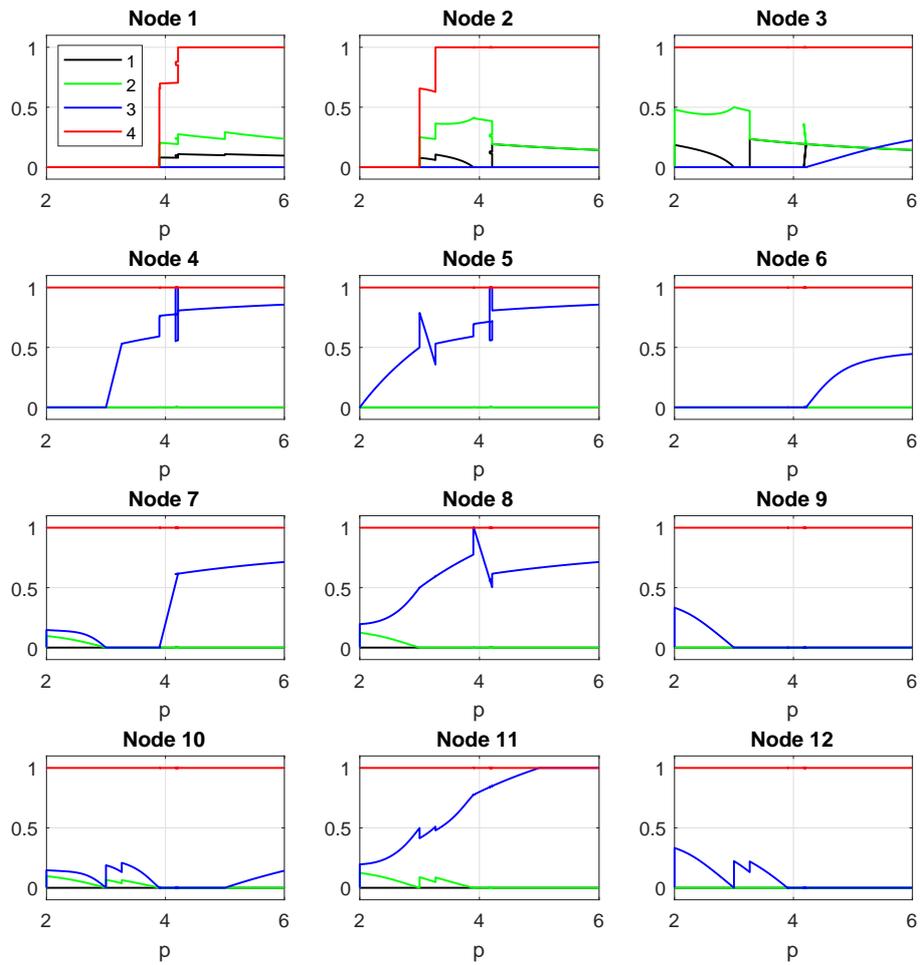}
 \caption{Computed equilibrium betting frequencies for $N=4$.}\label{freq_N4_2}
 \end{center}
 \end{figure}

Figure~\ref{E_N4_2} does indicate that $P=3$ (along with the other three values of $P$ shown) is a special case. There is also a small range of values of $P$ for which distinct multiple equilibrium solutions exist, although they have very similar expectations (almost indistinguishable in Figure~\ref{E_N4_2}). These can be seen more clearly close to $P = 4.1$ in Figure~\ref{freq_N4_2}, which shows the equilibrium betting frequencies as a function of $P$.
We can see from Figure~\ref{freq_N4_2} that for $P<3$, Players 1 and 2 check all their holdings at equilibrium, and the game revolves around Player 3 betting and the other players deciding whether to call with the second best card (Nodes 4, 5, and 6). Note that Nodes 7 to 12 are not reached at equilibrium for $P<3$. As $P$ increases past three, Player 2 starts to bet  at a nonzero frequency (Node 2), and for larger values of $P$, Player 1 becomes active at Node 1. The solution structure when $N=4$ is relatively simple because there is a clear distinction between bluffing cards (1 and 2), a single value card (4) and a single calling card (3). The only tension in the system is between value betting and sandbagging with the best card (4).

\subsection{Three player Kuhn poker with $N>4$ cards}
\begin{figure}
 \begin{center}
 \includegraphics[width=\textwidth]{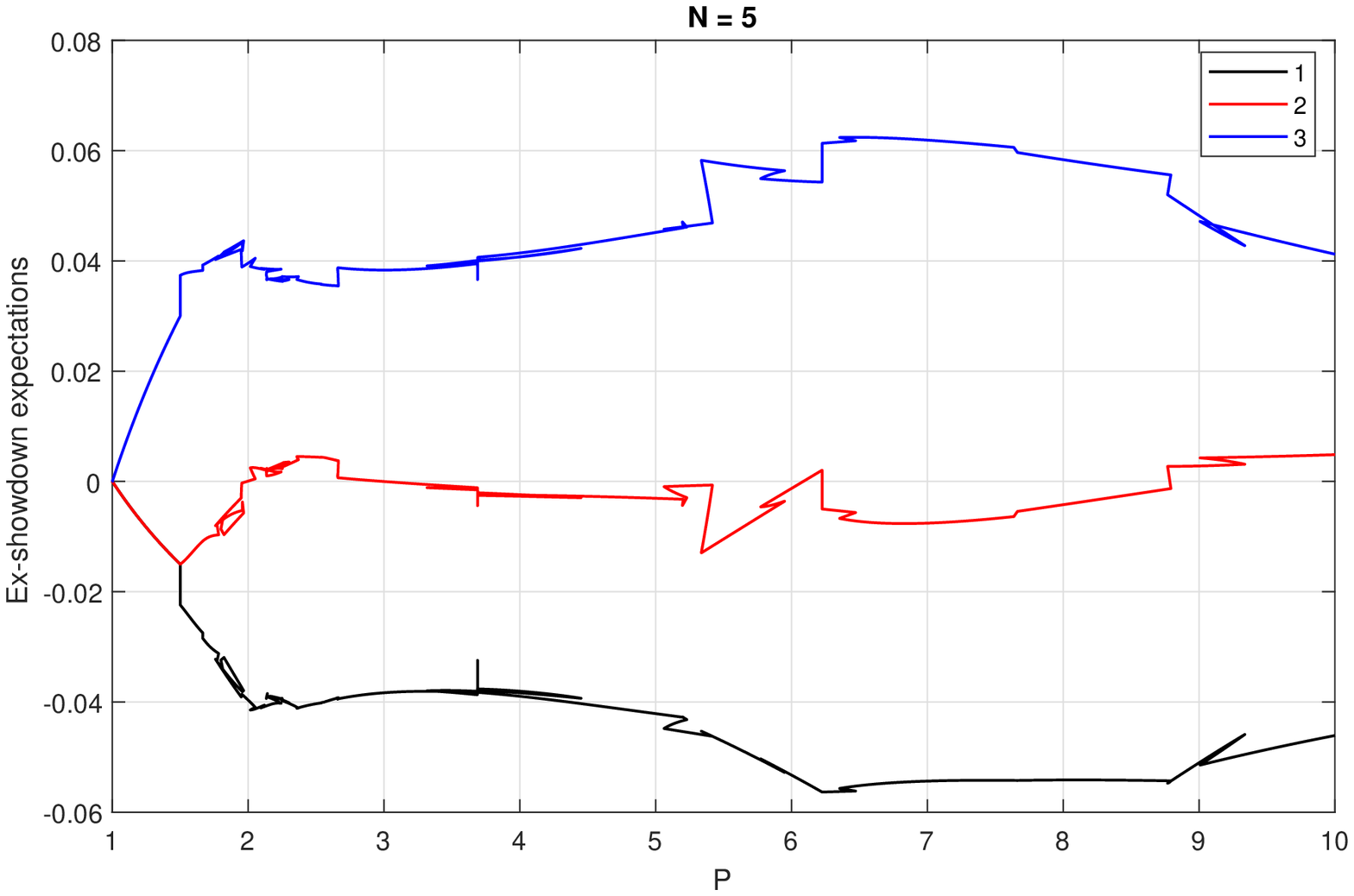}
 \caption{Computed expectations for $N=5$.}\label{E_N5}
 \end{center}
 \end{figure}

\begin{figure}
 \begin{center}
 \includegraphics[width=\textwidth]{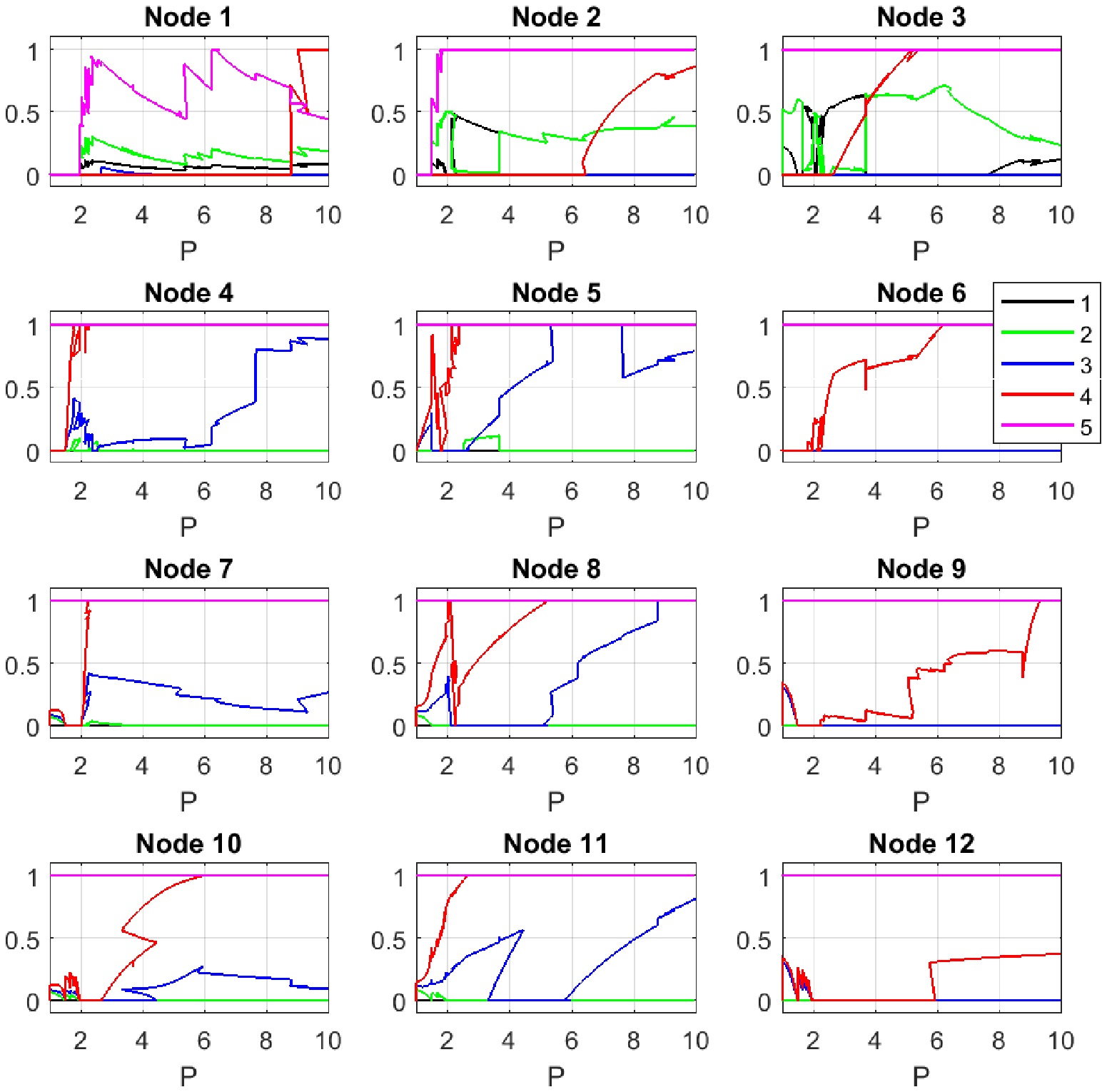}
 \caption{Computed equilibrium betting frequencies for $N=5$.}\label{freq_N5}
 \end{center}
 \end{figure}

Figure~\ref{E_N5} shows the equilibrium expectation when $N=5$. We can see that the addition of an extra card significantly increases the complexity of the structure of the solution curve (compare Figure~\ref{E_N5} with Figure~\ref{E_N4_2}). There are several ranges of values of $P$ for which there exist distinct, mutiple equilibrium solutions with a multiplicity of up to five, as opposed to three for $N=4$, and with differences in expected value that are more striking than those shown in Figure~\ref{E_N4_2}.
Figure~\ref{freq_N5} shows the equilibrium betting frequencies when $N=5$. Although there remains a clear distinction between value betting, bluffing and calling ranges, the behaviour at Node 1 is more complex than when $N=4$, with the value betting frequency of Player 1 with card 5 varying in an irregular manner for $P<10$, perhaps because now each player's value betting range may include the best two cards. When $N=4$ only the best card is ever value bet at equilibrium.

\begin{figure}
 \begin{center}
 \includegraphics[width=\textwidth]{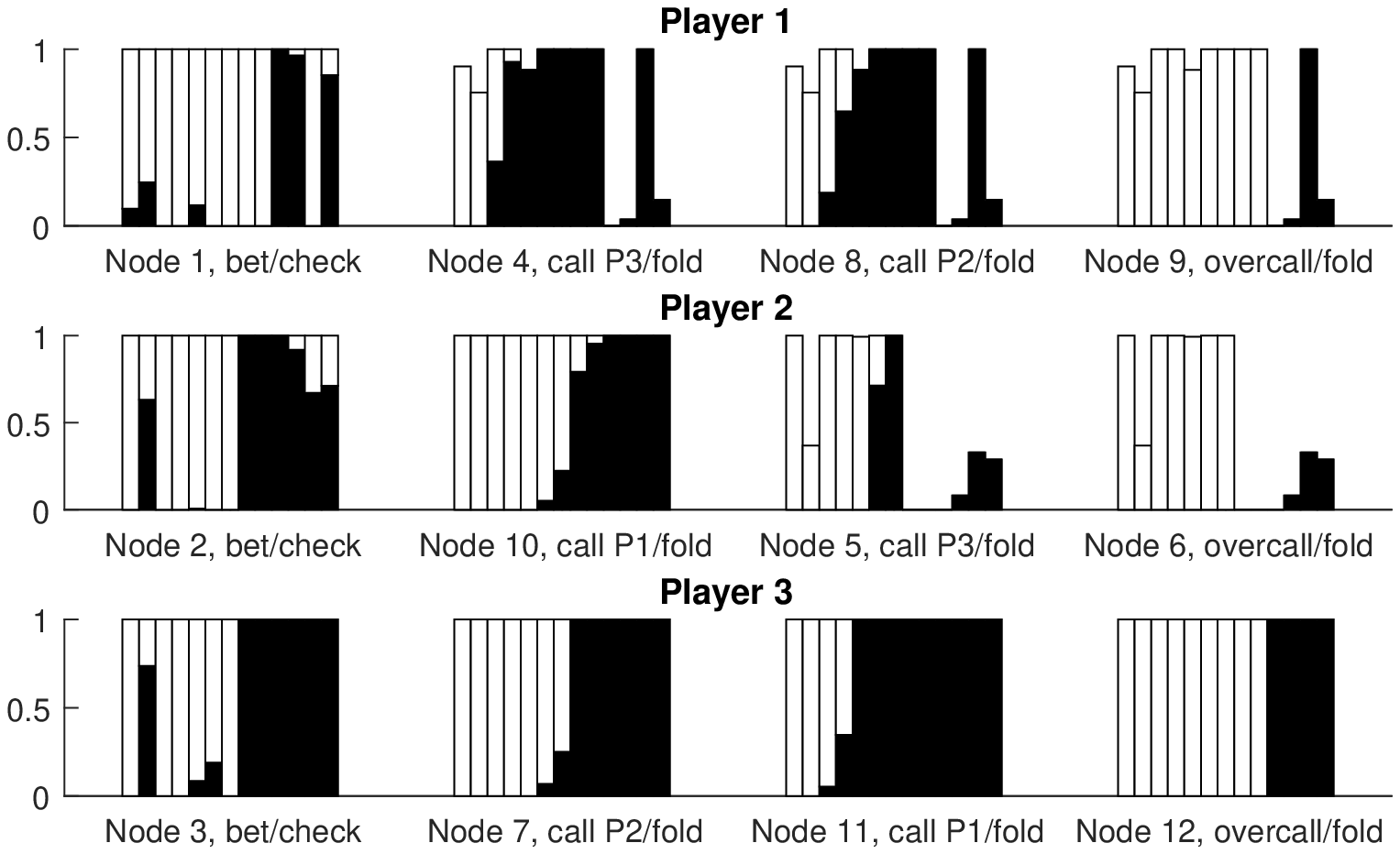}
 \caption{One of the computed sets of equilibrium betting frequencies for $N=13$ when $P = 9.2$. Each bar represents a card, ordered left to right, from lowest value to highest value. The shaded part of the bar indicates the frequency at which an aggressive action (betting or calling) is taken at equilibrium, whilst the height of the bar indicates the fraction of that holding that exists in each players equilibrium range at that node.}\label{freq_N13}
 \end{center}
 \end{figure}

For $N>5$ a new strategy emerges. The players, starting with just Player 1 when $N=6$, can, at equilibrium, make a bet with a mid-value card that is a bluff against one opponent and a value bet against  the other. This happens for a very small range of values of $P$ when $N=6$, so we will illustrate this for a game with more cards in the deck, $N=13$, where this phenomenon can be seen more clearly, as shown in Figure~\ref{freq_N13}. The case $N=13$ is also of interest because this is a game that can be played with a single suit of a real deck of cards.
 Figure~\ref{freq_N13} displays an equilibrium solution in a new format, for a single value of $P$. Each bar represents one of the cards, arranged from left to right with the highest value card at the right. The black part of each bar indicates the betting frequency at that node with that card. For example, in Figure~\ref{freq_N13}, at Node 2, which is controlled by Player 2, all of the best seven cards are bet with some nonzero frequency at equilibrium, balanced by some bluffs with the card of value 2. The bars at some of the subsequent nodes are shortened to indicate that some fraction has been bet at an earlier node. For example, at Node 5, Player 2 always calls with the three highest value cards, but the bars are shortened to indicate the fraction of those cards that are checked at equilibrium at Node 2. The cards with values 5, 6 and 7 do not appear in this equilibrium strategy at Node 5 since they are always bet by Player 2 at Node 2. 

In contrast to this clear bifurcation into value betting and bluffing ranges at Node 2, consider Node 3, controlled by Player 3, and the subsequent Nodes 4 and 5, where Players 1 and 2 must choose a calling strategy. Player 3 bets the cards of value 5 and 6 at a nonzero frequency at equilibrium, which forms a small but distinct range between the value betting and bluffing ranges. These bets are for value against Player 1 (Player 1 can call with worse at Node 4, but never folds a better card) and are bluffs against Player 2 (Player 2 never calls with worse at Nodes 5, but sometimes folds a better card). These two-way bets are a feature that can have no counterpart in two player games.

Although plots such as Figure~\ref{freq_N13} work well as a method of displaying a single equilibrium solution, they do not give any sense of how the equilibrium strategy varies along the solution curve. A better display method would be a succession of these plots. Since this would require far too many plots to include here, we have created videos for several different values of $N$ that illustrate the equilibrium solutions in the format of Figure~\ref{freq_N13}, with one frame for each solution, moving along the solution curve in equal steps of arc-length ($\left\| {\bf X}_i-{\bf X}_{i-1}\right\|$), linearly interpolating the numerically-calculated solution to achieve this. The videos are play at four frames per second by default, and are available online \footnote{The link to the videos is \url{https://drive.google.com/open?id=1cv1e-w4VqVd7Q56_Qa9comE01-KzBqol}, where CSV data files can also be found.}.

\begin{figure}
 \begin{center}
 \includegraphics[width=\textwidth]{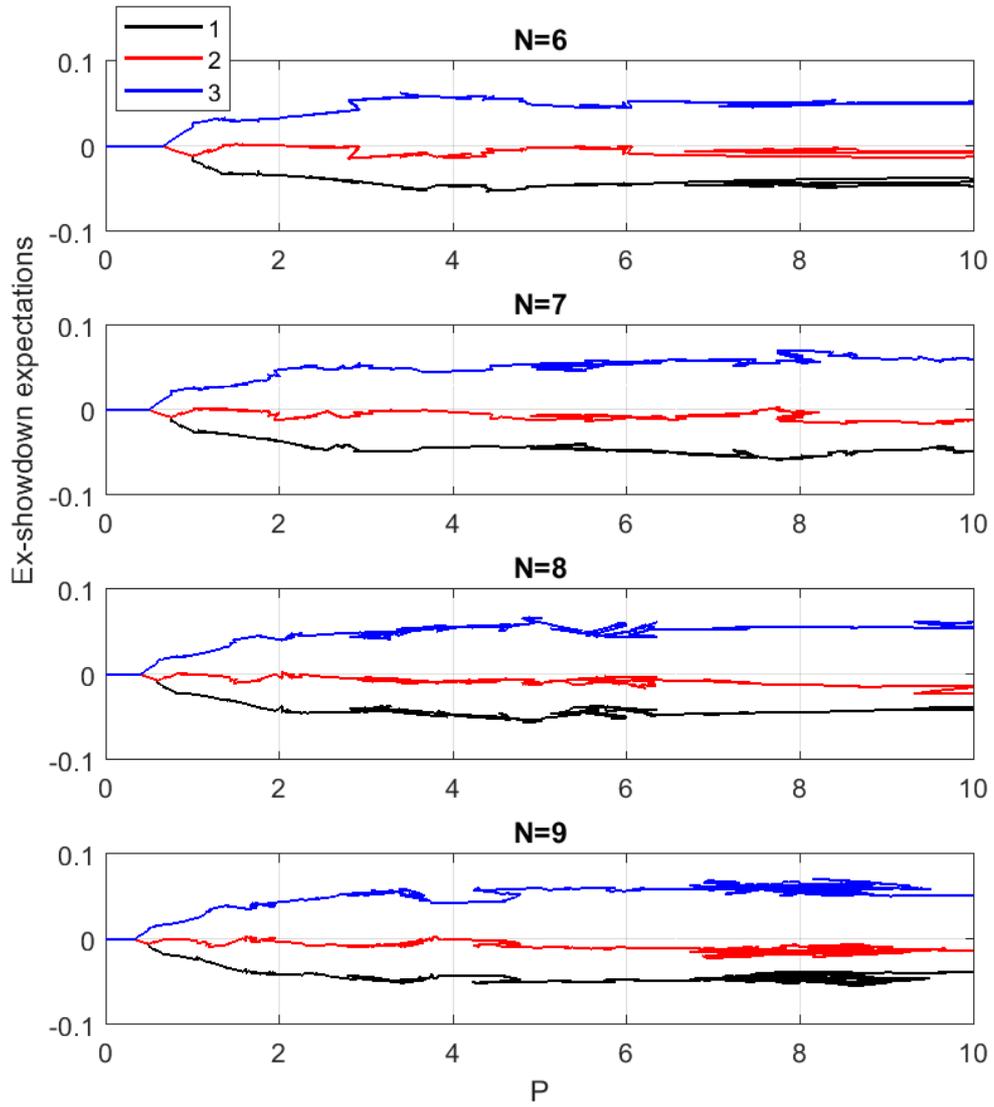}
 \caption{Computed expectations for $N=6$, 7, 8 and 9.}\label{E_N6-9}
 \end{center}
 \end{figure}

\begin{figure}
 \begin{center}
 \includegraphics[width=\textwidth]{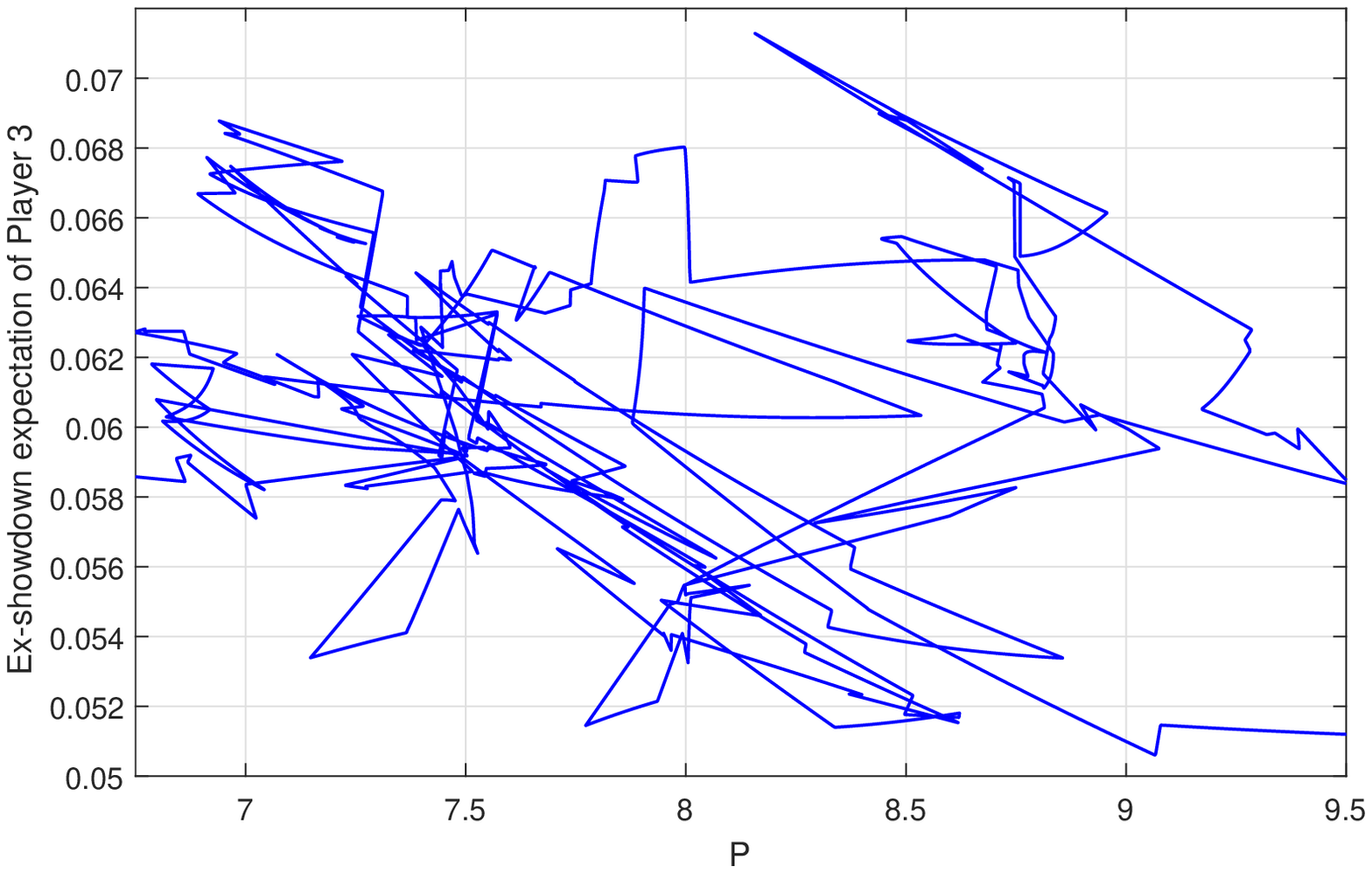}
 \caption{Computed expectation of Player 3 for $N=9$.}\label{E_N9}
 \end{center}
 \end{figure}

\begin{figure}
 \begin{center}
 \includegraphics[width=\textwidth]{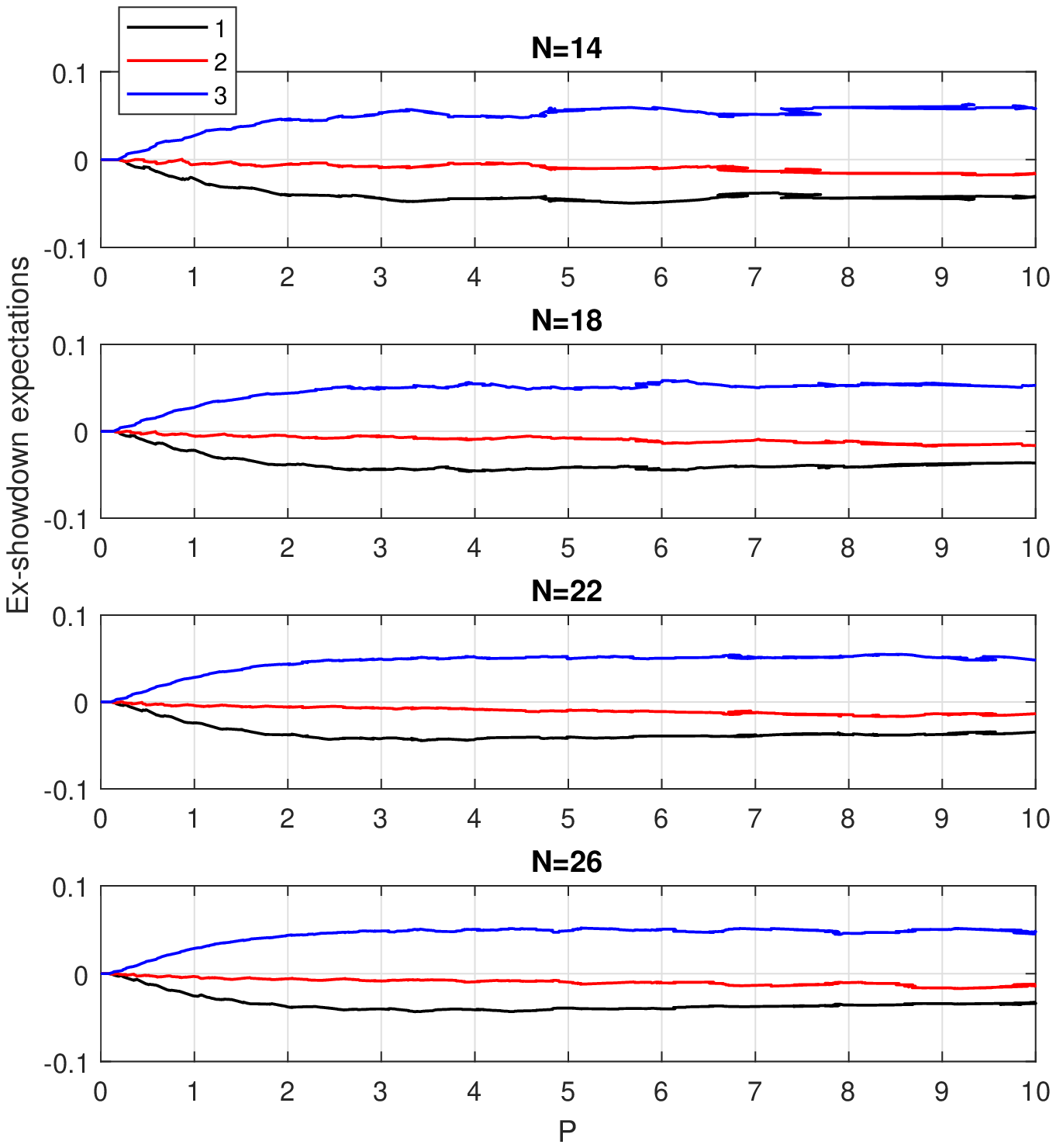}
 \caption{Computed expectations for $N=14$, 18, 22 and 26.}\label{E_N14-26}
 \end{center}
 \end{figure}

\begin{figure}
 \begin{center}
 \includegraphics[width=\textwidth]{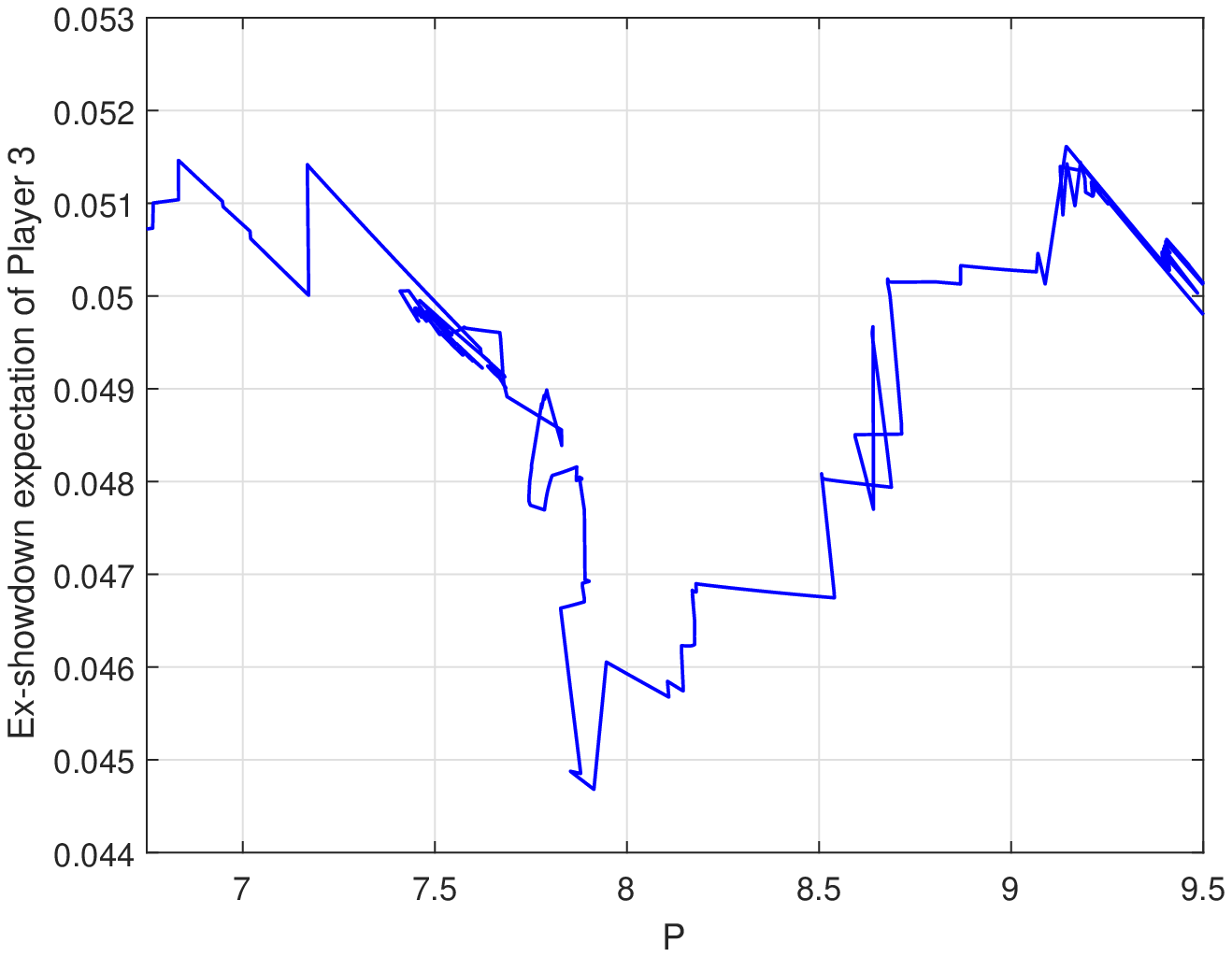}
 \caption{Computed expectation of Player 3 for $N=26$.}\label{E_N26}
 \end{center}
 \end{figure}

Finally, we can illustrate the complexity of the structure of the equilibrium solution curve by plotting the computed expectations for various values of $N$. Figure~\ref{E_N6-9} shows our results for $N = 6$ to $N=9$. The solution curves are hard to distinguish on this scale, so we have plotted Figure~\ref{E_N9}, which is a close up in the neighbourhood of $P=8$ when $N=9$, showing the expectation of Player 3. In the neighbourhood of $P=7.5$ there can be up to twenty four distinct, coexisting equilibrium solutions.
Figure~\ref{E_N14-26} shows the expectations for some larger values of $N$. Comparing this to Figure~\ref{E_N6-9} we can see that each curve has a smaller range of values at each value of $P$, which suggests that, as $N$ increases, although there is still a complex solution structure, each equilibrium solution for a given value of $P$ leads to a similar expectation. In Figure~\ref{E_N26} we can see that, if we zoom in on the expectation for $N=26$ (a game that could be played with two suits of a deck of cards), on this fine scale, the structure of the solution is still complex, despite the tightening of the range of values of the expectation.

The case $N=26$ is the largest game that we are able to compute with the memory that we have available (about 10Gb of RAM when $N=26$) so this method is not scaleable to even moderately large games. The bottleneck is the use of symbolic algebra to compute the expectations and their derivatives, which leads to very large MATLAB function files that use large amounts of memory when evaluated. After these functions have been computed, however, the method is very efficient. Data files containing the equilibrium solutions in CSV format are available online, along with the video files discussed earlier, for a range of values of $N$.

\section{Conclusions}\label{sec_conc}
In this paper, we introduced a new method to compute equilibrium solutions of extensive form, multiplayer toy poker games with no more than two branches at each decision node. After regularizing the game, we computed equilibrium solutions as a function of pot size, $P$,  using Newton's method and arc-length continuation. We used symbolic algebra in MATLAB to precompute the expectations and their derivatives. This leads to a numerical bottleneck in memory, as the symbolically computed functions are not memory efficient, and the largest game that we are able to deal with involves just twenty six cards. This problem can be overcome by computing the expectations and their derivatives directly from the game tree (with some loss of computational efficiency). The polynomial functions that define the expectaions, $E_i$, are linear in each betting frequency, so, rounding error notwithstanding, derivatives can be computed accurately using finite differences since perturbations of unit size lead to exact derivatives. For this to be done efficiently, we must use a compiled language and an algorithm that uses the structure of the game tree to minimize computation. This is currently under investigation. 

The numerical method discussed in this paper can be used equally well in the equivalent game with four or more players. We have investigated this using our current approach, but the memory requirements are very severe, and it is not feasible to compute solutions of even, for example, the four player game with six cards, so we will discuss these games in a later paper. Another extension to these games that would make them more relevant to standard limit poker is to allow up to three branches at each node (i.e. after a bet, a player is able to fold, call or raise). We have developed an extension of our regularization method that works for games with more than two branches per decision node, for which the mathematical formulation is more complicated than that given by (\ref{eqn1}), but again, we need to avoid using symbolic algebra, and this will form the subject of a later paper.

It is worth noting that, as discussed in \cite{SKP2017}, finding the equilibrium solutions of a three player game is an important first step in understanding the game, but, unlike the case of two player, zero sum games, simply using an equilibrium strategy (which equilibrium strategy if more than one exists?) does not guarantee unexploitability. Although the location of the equilibrium strategies likely indicates the approximate region of parameter space in which players should choose their strategies, the dynamic interaction of the players from hand to hand is likely to be chaotic, and dynamic strategies optimal.

Finally, an interesting technical question is whether the large $N$ limit of the games discussed here maps onto an equivalent $[0, 1]$ game. The equivalent $[0, 1]$ game is obvious for the half street, two player version of the game: Player 2 has continuous value betting and bluffing ranges, whilst Player 1 has a continuous calling range. Even for the two player game however, if there is a full street of betting and Player 1 must split the best cards into value betting and sandbagging ranges, the mapping to an equivalent  $[0, 1]$ game as $N \to \infty$ is not obvious, and merits further investigation.

\begin{appendix}
\section{Asymptotic solution for $0 < \epsilon \ll1$}\label{sec_asol}
In order to investigate the limit $\epsilon \to 0$ of the regularized problem (\ref{eqn2}), we construct the formal asymptotic solution. We order the equations in (\ref{eqn1_2}) so that an equilibrium solution, $x_k = \hat{x}_k$, satisfies
\begin{equation}
\begin{array}{ll} 
\hat{x}_k = 0 & \mbox{for $1 \leq k \leq k_0$,}\\
\hat{x}_k = 1 & \mbox{for $k_0+1 \leq k \leq k_0+k_1$,}\\
0 < \hat{x}_k < 1 & \mbox{for $k_0+k_1+1 \leq k \leq 12N$.}
\end{array}
\label{eqn_sol}
\end{equation}
There are thus $k_0$ equilibrium betting frequencies equal to zero, $k_1$ equal to one, and $k_2 = 12N - k_0 - k_1$ that lie strictly between zero and one. Note that, since (\ref{eqn_sol}) is an equilibrium solution,
\begin{equation}
f_k(\hat{\bf x})  \left\{
\begin{array}{ll} 
<0 & \mbox{for $1 \leq k \leq k_0$,}\\
>0 & \mbox{for $k_0+1 \leq k \leq k_0+k_1$,}\\
=0 & \mbox{for $k_0+k_1+1 \leq k \leq 12N$.}
\end{array}
\right.\label{eqn_expansion}
\end{equation}

We now expand the solution as
\begin{equation}
x_k = \hat{x}_k + \epsilon \bar{x}_k, \label{eqn_expansion2}
\end{equation}
with $\bar{x}_k = O(1)$ as $\epsilon \to 0$ and substitute into (\ref{eqn2_2}). At leading order as $\epsilon \to 0$, we immediately obtain the solutions
\begin{equation}
\bar{x}_k =  \left\{
\begin{array}{ll} 
-k_- /f_k(\hat{\bf x}) & \mbox{for $1 \leq k \leq k_0$,}\\
-k_+/ f_k(\hat{\bf x})  & \mbox{for $k_0+1 \leq k \leq k_0+k_1$.}
\end{array}
\right.\label{eqn_asol1}
\end{equation}
The remaining $k_2$ unknowns, $\bar{x}_k$, satisfy the $k_2$ linear equations
\begin{equation}
\sum_{l=k_0+k_1+1}^M f_{k,l}(\hat{\bf x}) \bar{x}_{l} = g^{-1}(\hat{x}_k)+k_-\sum_{l=1}^{k_0}\frac{ f_{k,l}(\hat{\bf x}) }{f_{l}(\hat{x})}+k_+\sum_{l=k_0+1}^{k_0+k_1}\frac{ f_{k,l}(\hat{\bf x}) }{f_{l}(\hat{x})} ~~\mbox{for $k_0+k_1+1 \leq k \leq 12N$,}\label{linsys}
\end{equation}
where $f_{k,l} \equiv \partial f_k/ \partial x_l$, noting that the function $g(y)$ is invertible since it is monotonically increasing. The asymptotic solution has this simple structure if and only if the function $g(y)$ satisfies (\ref{g_asym}), which is why this is a crucial requirement for this regularization method to be useful.

 If the matrix $f_{k,l}(\hat{\bf x}) $ in (\ref{linsys}) is non-singular, this determines $\bar{x}_l$ for $l = k_0+k_1+1, \ldots, 12N$. If $f_{k,l}(\hat{\bf x}) $  is singular, its nullspace corresponds to components of the equilibrium solution that are not uniquely-determined (i.e. they satisfy inequality constraints). In this case, the solvability condition provides equations that determine the leading order equilibrium solution. Since these equations are nonlinear, there is no obvious way of showing that they have a unique solution, but this is what we have found in all of our numerical computations. In other words, the regularization appears to select a unique equilibrium solution from those available whenever there is an inequality constraint in the exact solution.

We can illustrate this analytically using SKP, which has, in the notation of  \cite{SKP2017_2}, ordering the unknowns as ${\bf x} = \left(c_2, d_3, b_1, a_1, c_3, d_1, b_2, a_2, c_1, d_2, b_3\right)$,
\begin{equation}
\begin{array}{l}
f_1 = Pb_1-2a_1,\\
f_2 = (P+1)b_1-2a_1,\\
f_3=2P-4-(P+1)(c_2+d_3),\\
f_4=c_2+d_3-b_3-(1+\frac{1}{2}c_3)b_2,\\
f_5=(P+a_1)b_2+(b_1-2)a_2,\\
f_6=(P+1)b_2-2a_2,\\
f_7=2P-4+2a_1-(P+1)(c_3+d_1),\\
f_8 = c_3+d_1-\frac{1}{2}c_3b_1-(1+\frac{1}{2}c_1)b_3,\\
f_9=(P+a_2)b_3+b_2-2,\\
f_{10} = (P+1)b_3+b_1-2,\\
f_{11}=2P-4+2a_1+2a_2-(P+1)(c_1+d_2).
\end{array}
\end{equation}
Solution 1 exists for $2 \leq P \leq 3$ and is given by
\[
a_1 = b_1 =a_2=b_2=c_1= 0,~~b_3 = \frac{2}{P+1},~~d_2 = \frac{2P-4}{P+1},
\]
\begin{equation}
\frac{2P-4}{P+1} \leq c_2 + d_3 \leq \frac{2}{P+1},~~\frac{2P-4}{P+1} \leq c_3 + d_1 \leq \frac{2}{P+1}.
\end{equation}
Using hats for the leading order equilibrium solution and bars for the correction, from the first four equations we obtain
\begin{equation}
\begin{array}{l}
P\bar{b}_1 - 2 \bar{a}_1 = g^{-1} \left(\hat{c}_2\right),\\
(P+1)\bar{b}_1 - 2 \bar{a}_1 = g^{-1} \left(\hat{d}_3\right),\\
\bar{b}_1 = k_-\left\{(P+1)(\hat{c}_2+\hat{d}_3)-(2P-4)\right\}^{-1},\\
\bar{a}_1 = k_-\left\{\frac{2}{P+1}-(\hat{c}_2+\hat{d}_3)\right\}^{-1}.
\end{array}\label{a1b1eqns}
\end{equation}
The next four equations show that $\bar{a}_2$, $\bar{b}_2$, $\hat{c}_3$ and $\hat{d}_1$ also satisfy a system equivalent to (\ref{a1b1eqns}). Finally we have
\begin{equation}
\begin{array}{l}
\bar{c}_1 = \frac{1}{2}k_-(P+1),\\
(P+1)\bar{b}_3 = -\bar{b}_1+g^{-1}\left(\frac{2P-4}{P+1}\right),\\
(P+1) \bar{d}_2 = 2\bar{a}_1+2\bar{a}_2-(P+1)\bar{c}_1-g^{-1}\left(\frac{P}{P+1}\right),
\end{array}
\end{equation}
which give $\bar{c}_1$, $\bar{b}_3$ and $\bar{b}_2$ explicitly in terms of the other correction terms. By manipulating (\ref{a1b1eqns}), we find that
\begin{equation}
X = g\left(\frac{k_- \frac{P}{P+1}}{X-\frac{2P-4}{P+1}}-\frac{2k_-}{\frac{2}{P+1}-X}\right)+g\left(\frac{k_- }{X-\frac{2P-4}{P+1}}-\frac{2k_-}{\frac{2}{P+1}-X}\right) \equiv F(X),\label{eqn_X}
\end{equation}
where $X = \hat{c}_2+\hat{d}_3= \hat{c}_3+\hat{d}_1$. Since $dF/dX < 0$ and $g$ is monotonically increasing, (\ref{eqn_X}) has a unique solution for any $k_->0$ and $2 < P < 3$ by the Intermediate Value Theorem.

\end{appendix}

\newpage
\bibliography{Kuhn3}{}
\bibliographystyle{hplain}

\end{document}